\documentclass[12pt]{amsart}
\usepackage{amsmath,amssymb,amsthm,amscd,amsopn,amsfonts,amsxtra} 
\usepackage{latexsym,array}
\usepackage[mathscr]{eucal}
\usepackage{graphicx,psfrag}
\usepackage{epsfig}
\usepackage{color}
\newtheorem{theorem}[subsection]{Theorem} 
\newtheorem{lemma}[subsection]{Lemma}
\newtheorem{corollary}[subsection]{Corollary}
\newtheorem{proposition}[subsection]{Proposition}

\definecolor{orange}{rgb}{0.995, 0.75, 0.35}
\definecolor{purple}{rgb}{0.7, 0.2, 0.5}
\definecolor{royalblue}{rgb}{0.2, 0.7, 0.8}
\def\al{\alpha}
\def\de{\delta}
\def\eps{\epsilon}

\def\lam{\lambda}

\def\De{\Delta}

\def\th{\tanh}
\def\sech{\mathrm{sech}}
\def\supp{\mathrm{supp}}

\newcommand{\la}{\langle}
\newcommand{\ra}{\rangle}
\newcommand{\nd}{\noindent}

\newcommand{\vs}{\vspace}

\newcommand{\hB}{\hfill$\Box$}

\newcommand{\rk}{{\bf Remark.}\ \ }

\newcommand{\Z}{\mathbb{Z}}
\newcommand{\R}{\mathbb{R}}
\newcommand{\C}{\mathbb{C}}
\newcommand{\N}{\mathbb{N}}

\begin{document}
\title[Spectral multipliers for Schr\"odinger operators]
{Spectral multipliers for  Schr\"odinger operators with P\"oschl-Teller potential} 
\author{Shijun Zheng}
\address[]{Department of Mathematics \\
         Industrial Mathematics Institute\\
         University of South Carolina\\
         Columbia, SC 29208}
\email{shijun@math.sc.edu}
 \urladdr[]{http://www.math.sc.edu/\symbol{126}{shijun}}
\thanks{The author is supported by DARPA grant HM1582-05-2-0001} 

\keywords{spectral multiplier, Schr\"odinger operator, Littlewood-Paley theory}
\subjclass[2000]{Primary: 42B25; Secondary: 35J10, 35P25, 35Q40}
\date{September 29, 2006}

\begin{abstract}

We prove a sharp Mihlin-H\"ormander 
multiplier theorem for Schr\"odinger 
operators $H$ on $\R^n$. 
The method, which allows us to deal with general potentials, 
  improves Hebisch's method relying on heat kernel estimates for
positive potentials \cite{He90a,DOS02}.   
Our result applies to, in particular, the {\em negative} P\"oschl-Teller potential  
 $V(x)= -\nu(\nu+1)\,\sech^2 x $, $\nu\in \N$, for which $H$ has a resonance at zero. 
\end{abstract}

\maketitle 

\section{Introduction}\label{S1}

Spectral multiplier theorem for differential operators plays a significant role in 
harmonic analysis and PDEs. It is closely related to 
the study of the associated function spaces and Littlewood-Paley theory. 
Let $H=-\De+V $ be a Schr\"{o}dinger operator on $\R^n$, 
 where $V$ is real-valued. Spectral multipliers for $H$ have been considered in
\cite{He90a,E96,D99,D01,BZ05} and \cite{DOS02} for 
positive potentials.  
The case of negative potential is quite different and is not covered by the methods in
these papers.  Resonance and eigenvalue can occur that makes the analysis more involved.
 In this paper we are mainly concerned with proving a Mihlin-H\"ormander type
multiplier theorem on $L^p$ spaces for the Schr\"odinger operator with the negative P-T potential
\begin{equation}\label{e:p-t-nu}
V(x)= -\nu(\nu+1)\,\sech^2 x , \qquad \nu\in \N .
\end{equation}
In \cite{OOPSZ,Z06b} we are able to extend the sharp spectral multiplier theorem on 
 Triebel-Lizorkin 
spaces by modifying the argument in this note.  


Spectral multiplier problem requires both high and low energy analysis.
In high energy  the kernel of the multiplier operator $m(H)$ can be controlled by 
a weighted $L^2$ estimates.  In low energy, roughly, it can be controlled pointwise by 
an approximation to the identity.  This is the idea when dealing with positive selfadjoint
operators \cite{He90a,DOS02} where there is a 
rough kernel. 

However for negative operator in the Schr\"odinger case, resonance and eigenvalue 
can occur even for smooth and rapidly decaying potentials, which lead to failure of
the pointwise control of the kernel in lower energy.  The purpose of this paper is
to develop a general treatment to overcome this difficulty.  We find that this pointwise 
estimate can be substituted with a weaker estimate (in integral form) that turns out still work. 
This is the approach we will apply to the P\"oschl-Teller potentilal model.
The P-T potential arises in standing wave problem for the cubic wave equation 
\begin{equation}\label{eq:nlw}
-u_{tt}+u_{xx}+2u-u^3=0 \, . 
\end{equation}

In  \cite{Z06b} we give a general treatment on spectral multiplier
problem for Schr\"odinger operators satisfying for every $j\in \Z$
\begin{equation}\label{e:phi-dec}
\vert  \Phi_j( H)(x,y)\vert \le c_n\frac{2^{nj/2}}{(1+2^{j/2}|x-y|)^{n+\eps}}.
\end{equation}
where $\Phi_j(x)=\Phi(2^{-j}x)$, $\Phi\in C^\infty_0(\R)$. 

The assumption is verified when $H$ is a Schr\"odinger operator $-\De+V$,  
$V\ge 0$ is in $L^1_{loc}(\R^n)$ \cite{Ou06, He90a} or
$H$ is a uniformly elliptic operator on $L^2(\R^n)$ 
\cite[Theorem 3.4.10]{Da89}. It is showed in \cite{He90a,Z06a} 
that the decay  (\ref{e:phi-dec}) is satisfied whenever 
the heat kernel of $e^{-tH}$ satisfies the upper Gaussian bound
\begin{equation*}
 0\le e^{-tH}(x,y) \le c_n t^{-n/2} e^{-c |x-y|^2/t} \,. 
\end{equation*}
However, when $V$ is negative, eigenvalue and resonance may occur at the origin. 
the seemingly ubiquitous decay (\ref{e:phi-dec}) for general selfadjoint operators is 
not valid for all $j$. Our approach shows that if (\ref{e:phi-dec}) replaced
by an integral version (\ref{e:int-phi-dec}), the argument in \cite{He90a, DOS02} can still work.   
This treatment will be further elaborated in the study of spectral calculus for 
rough potentials in the critical class in a following paper.  

 The basic ingredients we need to show are two weighted inequalities
\begin{itemize}
\item[(a)] \label{e:int-phi-dec} 
If $\Phi\in C^\infty_0(\R)$ there exists a finite measure $\mu$ such that for each interval $I$ with length $2^{-j/2}$, 
  $j\in\Z$,
\begin{align*}  |\Phi_j(H)(x,y)|
 \le& c \int_{u\in\R^n} \frac{2^{jn/2}}{(1+2^{j/2}| x-y-u|)^{n+\eps}} d\mu(u)\\
\approx&  \rho_j *\mu \, (x-y)  \notag\end{align*} 
 where 
$\rho_j(x):= 2^{jn/2} (1+ 2^{j/2}|x|)^{-n-\eps}$.
\item[(b)] 
\[
\sup_y \Vert \la 2^j(x-y)\ra^\al \phi_j(H)(x,y)\Vert_{L^2_x}\le C  
\sup_j\Vert \chi \phi(2^{-j}\xi^2)\Vert_{X^\al}
\]
where $X^\al$ is $C^\al(\R)$ or $W^{\al}$, $\al>n/2$.  $\chi$ is a given smooth cut-off function
with support off $0$.
\end{itemize}

Observe that (a) is a weaker condition than the pointwise estimate for the decay of $\phi_j(H)$.
Indeed, the ponitwise decay (\ref{e:phi-dec}) is a special case when $\mu=\de$.

\subsection{Weighted $L^2$ estimates for kernel of $m(H)$}

 There are a few ways to prove $L^p$ boundedness for $m(H)$.




The usual condition is the H\"ormander integral condition
\begin{equation*}
 \int_{|x-y|>2|x-\bar{y}|} |m(H)(x,y)-m(H)(x,\bar{y})|dx\le A 
\end{equation*}
for $y\in I $, $\bar{y}$ the center of $I$, $I$ being any cube in $\R^n$. 
 This is what is shown in \cite{BZ05}. 
For P-T potential similar estimate is not valid 
in low energy.  
This is the reason why we need
to consider the weighted $L^2$ estimate in (b).  

The other way is to use wave operator method \cite{W99,Y95}.  However
wave operator method does not give the sharper weak $(1,1)$ result. 






\section{Main results}

Let $H$ be a selfadjoint operator on $L^2(\R^n)$. 
Then if $\phi\in L^\infty$, we can define
$\phi(H)=\int \phi(\lam) dE_\lam$ by functional calculus, where
$H=\int \lam dE_\lam$ is the spectral resolution of $H$.

Our main result is the following 

\begin{theorem}\label{th:m(H)Lp} 
Suppose  $H$ verifies the weighted decay (a) 
and weighted $L^2$ inequality (b).
Then $m(H)$ is bounded on $L^p(\R^n)$, $1<p<\infty$ and of weak type $(1,1)$.
Moreover, 
\[  \Vert m(H)\Vert_{L^1\to\text{w-}L^1} \le C(m) 
\]
$C(m):= \Vert m\Vert_\infty+\sup_{\lam>0}\Vert \chi(\cdot)m(\lam\cdot)\Vert_{X^\al }$\,.
\end{theorem}


\nd
\rk One may view condition (a) as a ``pointwise" control of the kernel in {\em lower} energy
while condition (b) as a norm control in {\em higher} energy. 


  The conditions of Theorem \ref{th:m(H)Lp} applies to 
the case when the kernel of $\phi(H)$ is slowly decaying or lack of smoothness.  
It can also simplifies the proof of some known results on other potential of
polynomial growth, e.g. the Hermite operator on $\R^n$ \cite{Z06b}. 

Applying the theorem to P-T potential we obtain

\begin{theorem}\label{th:p-t-m(H)} 
Let  $H=-d^2/dx^2 +V_\nu$, $\nu\in\N$,  where $V_\nu$ is the P-T potential in (\ref{e:p-t-nu}).
Then $H$ satifies the weighted decay (a) 
and weighted $L^2$ inequality (b) with $X^\al=C^\al$, $\al=1$. 
Therefore  the conclusion of Theorem \ref{th:m(H)Lp} holds for the 
 one dimensional P-T model. 
\end{theorem}

From section \ref{S2} we know the derivative of the kernel fail to satisfy nice decay, 
making it difficult to control the difference of $m(H)$ and thus
leading to the failing of low energy estimates for the H\"ormander integral condition.

A recently developed approach \cite{DOS02} by Sikora et al extended Hebisch's method \cite{He90a}
that apply to positive operators efficiently 
but rely heavily on heat kernel estimates. 
However, in dimensions one
and two when the potential is negative, such a heat kernel estimate is NOT 
available. Therefore we consider more direct approach and would rather state and 
prove the multiplier result in Theorem \ref{th:m(H)Lp} for general dimensions. 


We will assume $\Phi, \varphi\in C_0^\infty ({\R}) $ satisfy 
the condition 
\begin{align}
&\sum_{j=-\infty}^\infty \varphi_j(x) =1    \quad x\neq 0     \label{eq:dyadic-id}\\
&\Phi(x)  +\sum_{j=1}^\infty \varphi_j(x) =1,    \quad   \forall x
\end{align}
where $\varphi_j(x) =\varphi( 2^{-j} x)$. 

\section{Proof of Theorem \ref{th:m(H)Lp}} 

The technical lemmas 
we need in proving  Theorem \ref{th:m(H)Lp} are:
\begin{lemma} \label{l:j-horm-cond}
  Let $y\in I$, $I\subset\R^n$ a cube with length $t=\ell(I)$.  Let $2^{-j_I/2}\sim  t$.
Then
(a) For $t>0, j\in\Z$, 
$$\int_{|x-y|\ge 2t} | m_j(H)(1-\Phi_{j_I}(H))(x,y)| dx \le C(2^{j/2}t)^{\frac{n}{2}-s} 
\Vert m(\xi^2)\Vert_{\dot{H}^s}$$
$s>n/2$. 
(b)\[
\int_{|x-y|\ge 2t} \sum_{j=-\infty}^\infty| m_j(H)(1-\Phi_{j_I}(H) )(x,y)| dx \le A\,.
\]
In particular, 
\[
\int_{|x-y|\ge 2t} \sup_{j\in\Z} | m_j(H)(1-\Phi_{j_I}(H) )(x,y)| dx \le A\,.
\]
\end{lemma}





\begin{lemma}\label{l:max-min} Condition (b) implies 
\begin{align*}\max_{y\in I}  |\Phi_j(H)(x,y)|
 \le& c \min_{y\in I}\int_{u\in\R^n} \frac{2^{jn/2}}{(1+2^{j/2}| x-y-u|)^{n+\eps}} d\mu(u)\\
\approx& \min_{y\in I} \rho_j *\mu \, (x-y)  \notag\\
\le&   \frac{1}{|I|}\int_{z\in I} \rho_j *\mu \, (x-z)dz
\end{align*} 
 where 
$\rho_j(x):= 2^{jn/2} (1+ 2^{j/2}|x|)^{-n-\eps}$.\end{lemma}
 
 The proof is based on the observation if
  $\ell(I)=r_I$, the side length of a cube $I$ in $\R^n$, 
\[   \sup_{y\in I} (1+ |x-y|)^{-n-\eps}\le C \min_{y\in I} (1+ |x-y|)^{-n-\eps}\]
hence  
\[  \sup_{y\in I} (1+ |x-y|/t)^{-n-\eps}\le C \min_{y\in I} (1+ |x-y|/t)^{-n-\eps}\le
\frac{C}{|I|}\int_I  (1+ |x-y|/t)^{-n-\eps} dy\]
where $t\sim \ell(I)$, $I$ is any cube,   
see \cite{He90a}. 






\vs{.23in}
\nd
{\bf Proof of weak (1,1).}  If applying the C-Z decomposition we know the main part is how 
to handle the ``bad'' function $b=\sum_k b_k$
where $b_k\subset I_k$, $I_k$ being disjoint intervals in $\R$. 

\begin{proof} Let $\Phi\subset [-1,1]$, $\Phi_j(x)=\Phi(2^{-j}x)$.
Write 
\[ m(H)b(x)=\sum_k m(H)(1-\Phi_{j_k}(H)) b_{k}(x)
+ \sum_k m(H)\Phi_{j_k}(H) b_{k}(x).\] 
where $2^{-j_k}\sim \ell(I_k)^2$. 
We need to show
\begin{align*}
\vert&\{x\in\R\setminus \cup_kI_k^*: 
 |m(H) b(x)| >\lam/2 \}\vert\\
\le&  |\{x\in\R\setminus I_k^*:  \sum_k |m(H)(1-\Phi_{j_k}(H)) b_{k}(x)|>\lam/4\}\\
+&    |\{x\in\R\setminus I_k^*:  \sum_k |m(H)\Phi_{j_k}(H) b_{k}(x)|>\lam/4\} \\
\le &\lam^{-1} \Vert f\Vert_1,
\end{align*}
where $b=\sum_k b_k$ (convergence in $L^1\cap L^q$ so 
$T b(x)=  \sum_k Tb_{k}$ in $L^q$)

{\em Higher} energy 

Denote $I_k^*$ the cube having length three times the length of $I_k$ with the same
center as $I_k$.
 If $x\notin \cup_kI_k^*$, $I_k\subset \{y: |y-x|>r_k \}$, $r_k$ being the length of $I_k$.
\begin{equation*}
 m(H)(1-\Phi_{j_k}(H )b_k(x)=\int_{|y-x|>r_k} m(H)(1-\Phi_{j_k}(H)(x,y) b_k(y)dy
\end{equation*}
 Apply weighted condition Lemma \ref{l:j-horm-cond} (c)
\begin{align*}  
&|\{x\notin \cup I_k^*:  |\sum_k m(H)(1-\Phi_{j_k}(H )b_k(x)|>\al/4 \}| \\
 \le& C (\al/4)^{-1} \int_{\R^n\setminus\cup I_k^*}  |\sum_k m(H)(1-\Phi_{j_k}(H) )b_k(x)| dx \\
\le& C\al^{-1} \int  \sum_k |b_k(y)| dy\int_{|y-x|>r_k} |m(H)(1-\Phi_{j_k}(H)(x,y)|dx\\
 \le& C \sup_\lam\Vert \chi m(\lam\xi^2)\Vert_{C^s_{loc}} \al^{-1} \int |b(y)|dy \\
 \le& C \sup_\lam\Vert \chi m(\lam\xi^2)\Vert_{C^s_{loc}}\al^{-1}\Vert  f\Vert_1.
\end{align*}
where we note
\begin{align*}
&\int_{|x-y|>r_k} |m(H)(1-\Phi_{j_k}(H) )(x,y)| dx\\
\le& \sum_{2^j> r_k^2}\int_{|x-y|>r_k} |m_j(H)(x,y)| dx\le C
\end{align*}
because  if $\Phi(x)+\sum_{j=1}^\infty \phi(2^{-j}x)=1$
then for any  $j_0\in\Z$,
$\Phi(2^{-j_0}x)=1-\sum_{j=j_0+1}^\infty \phi(2^{-j}x)$.

{\em Lower} energy 

Since $m(H)$ is bounded on $L^2$. The proof is complete if we can show
\begin{equation}\label{e:sum-L2}
\int |\sum_k   \Phi_{j_k}(H)b_k(x)|^2 dx\le C \al  \Vert f\Vert_1 
\end{equation}
To show this let $h\in L^2(\R^n)$,  $2^{-j_k}\sim \ell(I_k)^2$. 
According to (b)   
\begin{align*} 
\la&  \sum_k \Phi_{j_k}(H)  b_k  , h \ra\\
=& \sum_k \int_x h(x)dx \int_{y\in I_k} K_{j_k}(x,y)  b_k(y)dy   \\
\le & \sum_k   
|I_k |^{-1}\int_x |h(x)|dx\int_{z\in I_k}\int_u  
\rho_j(x-z-u) d\mu(u) dz \int_y|b_k(y)| dy\\   
\le& 
\sum_k \Vert b_k\Vert_1 |I_k|^{-1}  \int_{z\in I_k}   \int_u (M h) (z+u)    d\mu dz  \\ 
&(bec \; \rho_j=2^{jn/2}(1+2^{j/2}(\cdot) )^{-n-\eps}  \; \text{is $\sim$ an approximation to the id})\\
\le& C \al \int_z \sum_k\chi_{I_k}(z)  (M h)* d\mu (z) dz\\
\le & C \al  \Vert\sum_k \chi_{I_k}\Vert_2 \Vert M h* d\mu \Vert_2\\
\le & C \al  (\sum_k |I_k|)^{1/2} \Vert h\Vert_2\\
\le & C \al  (\al^{-1} \Vert f\Vert_1)^{1/2} \Vert h\Vert_2
= C  \al^{1/2}\Vert f\Vert_1^{1/2} \Vert h\Vert_2
\end{align*}
which proves (\ref{e:sum-L2}).
We have used the fact that if  $\rho_t= t^{-n}\rho(x/t)$ is any approximation kernel to the identity,
so that $\rho\in L^1(\R^n)$ is positive and decreasing,  then 
\[
 \sup_{t>0} |\rho_t*f(x)| \le Mf(x)
 \]
where $M$ denotes the Hardy-Littlewood maximal function on $\R^n$.   
\end{proof}

\nd
\rk Note that in analyzing the kernel of $\Phi_j(H)$, oscillatory integral, if $j<0$  we can prove the rapid decay from the spectrum side; if $j>0$ it only gives a decay of $|x-y|^{-n}$; we will have to 
work in the space side and use the (average) integral version of rapid decay.  

\nd
\rk From the proof above we see the 
the weighted $L^2$ inequality in (b) somehow plays the role of   
H\"ormander condition in classical case \cite{St93,BZ05} 
\begin{equation}\label{e:int-hor-con}
\int_{|x-y|>2|x-\bar{y}|} |K(x,y)-K(x,\bar{y})|dx\le A 
\end{equation}
which requires the gradient estimate for $K_j(x,y)$.  
\nd
\rk A nontrivial vector-valued version of the  proof of the $L^p$ result 
 yields the multiplier result on the homogenous spaces
 $\dot{F}(H)$ and $\dot{B}(H)$ \cite{OOPSZ, Z06b}.

   In \cite{OZ} when identifying the inhomogeneous space
$F^{0,2}_p(H)=L^p$ we used decay for the derivative of the kernel 
in high energy, namely,  
(b') For $t>0, j\ge 0$,
$$\int_{|x-y|\ge 2t} | m_j(H)(x,y)-m_j(H)(x,\bar{y})| dx \le C(2^{j/2}t)^{\frac{1}{2}} $$
(c')\[
\int_{|x-y|\ge 2t} \sum_{2^j\ge t^2}^\infty |m_j(H)(x,y)| dx \le C.\]
which is valid only for high energy estimate because we only 
have available for $j\ge 0$ the weighted estimate 
\[ \Vert (x-y)\partial_y K_j(x,y)\Vert_2\le  2^{j/4}\,.
\]

To deal with the problem in low energy we avoid using the estimate for $\partial_yK_j(x,y)$
and follow the line of proof of Theorem 2.1 \cite{OOPSZ} 
for the homogeneous spaces 
$\dot{F}(H)$. Thus we obtain
the identification of $\dot{F}(H)$ and $L^p$ spaces.  
\begin{corollary}\label{c:F-Lp} Let $1<p<\infty$.  Then 
\[ \dot{F}_p^{0,2}(H)= \dot{F}_p^{0,2}(\R^n)=L^p(\R^n) ,\]
meaning that the Littlewood-Paley characterization holds
\[ \Vert\big(\sum_{j\in\Z} \vert\phi_j(H)f\vert^2\big)^{1/2}\Vert_p  \approx \Vert f\Vert_p\; .
\]
\end{corollary}

\section{Proof of Corollary \ref{c:F-Lp} }

\nd
  Identification of $\dot{F}^{\al,q}_p(H)=L^p$.   {\em homogeneous} spaces 

Let $Q_j=\phi_j(H)$, $j\in \Z$.  Define
\[ Q:  f\mapsto  \{ \phi_j(H)f  \}
\]
and 
\[  R:  \{f_j\}\mapsto \sum_{j=-\infty}^\infty \psi_j(H)f_j 
\]
We show that the same method in showing spectral multiplier theorem 
for $F$ spaces yields 
\[ Q:  L^1\to w-L^1(\ell^2)
\]
and 
\[ R:  L^1(\ell^2)\to w-L^1
\]
Hence, this, together with the boundedness $Q$: $L^2\to L^2(\ell^2)$ and 
$R:\; L^2(\ell^2)\to L^2$, proves that if $1<p<\infty$
\[     \Vert f\Vert_{L^p}\sim \Vert \phi_j(H)f\Vert_{ L^p(\ell^2)} = 
\Vert f\Vert_{ \dot{F}^{0,2}_p(H)} \]

\begin{lemma}\label{l:Q-R-L2} $Q$: $L^2(\R^n)\to L^2(\ell^2)$. 
$R:\; L^2(\ell^2)\to L^2(\R^n)$.
\end{lemma}

Proof. \begin{align*}  &\sum_{j\in\Z}(\phi_j(H)f,\phi_j(H)f )\\
=&(\sum_{j\in\Z}(\overline{\phi_j}\phi_j)(H)f,f )\lesssim \Vert f\Vert_2^2
\end{align*}
because $\sum_{j\in\Z}|\phi_j(x)|^2\approx 1$. 
The proof for R is similar. \hB

\begin{lemma}\label{l:L1-w-L1l2} Q is $L^1\to w-L^1(\ell^2)$, i.e.,
\[  |\{ \sum_{j\in\Z}   |\phi_j(H)f(x)|^2)^{1/2}>\al \}|\le C \al^{-1}\Vert f\Vert_1
\]
\end{lemma}

\subsection{Proof of $Q$: weak $(L^1,L^1(\ell^2))$}\label{ss:w-11}   Let $f\in L^1(\R^n) $.  
Given $\al>0$, let $f =g+b $  be the Calder\'on-Zygmund decomposition,
where $g\in L^2\cap L^1$, $b=\sum_k b_k$, $b_k\subset I_k$, $I_k$ disjoint.
\begin{itemize}
\item[i)] $\displaystyle{ |g(x)|\le C\al}$
\item[ii)] $\displaystyle{|I_k|^{-1}\int_{I_k} |f| dx\le C_n\al}$.
\end{itemize}
Note that (i), (ii) imply
\[  |I_k|^{-1}\int_{I_k} |b| dx\le C_n\al.
\]

Since $\phi_j(H)$  is bounded  $L^2\to L^2(\ell^2)$. 
\[  \int \sum_{j\in\Z} |\phi_j(H)g(x)|^2 dx \le C \Vert g\Vert_2^2\le C \al\Vert f\Vert_1
\] 
  (i)   gives 
\[  
|\{x: (\sum_{j\in\Z} |\phi_j(H)g(x)|^2)^{1/2}>\al/2     \}| \le 
C \Vert \phi\Vert_\infty\al^{-1}\Vert  f\Vert_1.
\]

We estimate $b$ in more detail.  Let $2^{-j_k}\sim t_k=$ diameter of the cube $I_k$. 
\begin{align*}  
&|\{x:  (\sum_{j\in\Z} |\phi_j(H)b(x)|^2)^{1/2}>\al/2     \}| \\
\le&|\{x: (\sum_j |\sum_k \phi_j(H)(1-\Phi_{j_k}(H)b_k(x)|^2)^{1/2}>\al/4     \}| \\
+& 
|\{x: (\sum_j |\sum_k  \phi_j(H)\Phi_{j_k}(H) b_k(x) |^2 )^{1/2}>\al/4     \}|\\
:=&B_1+B_2.
\end{align*}

For the first term since $|\cup I_k|\le   \al^{-1}\Vert  f\Vert_1$, it suffices to estimate: 

Denote $I_k^*$ the cube having length three times the length of $I_k$ with the same
center as $I_k$.
 If $x\notin \cup_kI_k^*$, $I_k\subset \{y: |y-x|>r_k \}$, $r_k$ being the length of $I_k$.
\begin{equation*}
 \phi_j(H)(1-\Phi_{j_k}(H)b_k(x)=\int_{|y-x|>r_k} 
 \phi_j(H)(1-\Phi_{j_k}(H)(x,y) b_k(y)dy
\end{equation*}
 Apply Hebisch-Zheng condition (b)
\begin{align*}  
&|\{x\notin \cup I_k^*:  (\sum_j |\sum_k  \phi_j(H)(1-\Phi_{j_k}(H) b_k(x)|^2)^{1/2}>\al/4 \}| \\
 \le& C (\al/4)^{-1} \int_{\R^n\setminus\cup I_k^*}  (\sum_j |\sum_k
 \phi_j(H)(1-\Phi_{j_k}(H) b_k(x)|^2)^{1/2} dx \\
\le& C\al^{-1} \int_{\R^n\setminus\cup I_k^*}
  \sum_k  (\sum_j   |\phi_j(H)(1-\Phi_{j_k}(H)b_k(x)|^2)^{1/2} dx\quad (Minkowski)\\
 \le& C  \al^{-1}  \sum_k\int |b_k(y)|dy 
 \int_{\R^n\setminus\cup I_k^*}\sum_{j\in\Z}|\phi_j(H)(1-\Phi_{j_k}(H)(x,y)|  dx \quad (Jensen)\\
\le& C\al^{-1}   \sum_k\int |b_k(y)|dy\int_{\R^n\setminus I_k}\sum_{j\ge j_k}|\phi_j(H)(1-\Phi_{j_k}(H)(x,y)|  dx\\
  \le& C \sup_\lam\Vert \phi_j(\lam\xi^2)\Vert_{H^s_{loc}}\al^{-1}\Vert  f\Vert_1.
\end{align*}


It remains to deal with $B_2$.  (Low energy)   
 The proof is finished if we can show
\begin{equation}\label{e:lower-b}
\int \sum_j | \phi_j(H)\sum_k \Phi_{j_k}(H) b_k(x)|^2 dx\le C \al  \Vert f\Vert_1 
\end{equation}
Since $Q=\{\phi_j(H)\}_{j\in\Z}$  
is bounded $L^2\to L^2(\ell^2)$, we only need to prove
\begin{equation}\label{e:sum-b-lower}
\int \sum_k |\Phi_{j_k}(H) b_k(x)|^2 dx\le C \al  \Vert f\Vert_1 \,.
\end{equation}
Similar to the proof of Theorem \ref{th:m(H)Lp} we use duality argument. 
Let $h\in L^2(\R^n)$, $\Vert h\Vert_{L^2}\le 1$. Then by condition (b) 
\begin{align*}
&<  \sum_k \Phi_{j_k}(H)  b_k(x)  , h>\\
\le&c\sum_k \Vert b_k\Vert_1 |I_k|^{-1}\int_{z\in I_k}dz\int_u d\mu(u)\int_x  \rho_{t}(x-z-u)  |h(x)| dx\\ 
=& c \sum_k |I_k|^{-1} \Vert b_k\Vert_1 
\int_z \chi_{I_k}(z)dz  \int_u d\mu(u)\int_x\rho_{t}(x-z-u)  |h(x)| dx\\
\le& c  \al \int_z \sum_k\chi_{I_k}(z) dz \int M h(z+u) d\mu(u)\\
\le & c \al\int_z \Vert\sum_k \chi_{I_k}\Vert_2 \Vert (M h)*\mu\Vert_2\\
\le & c  \al ( \al^{-1}\Vert f\Vert_1^2)^{1/2} \Vert h\Vert_2
\le c  (\al\Vert f\Vert_1)^{1/2} 
\end{align*}
which proves (\ref{e:sum-b-lower}) hence (\ref{e:lower-b}).  
This completes the proof of  Lemma \ref{l:L1-w-L1l2}. \hB

Similarly we can show

\begin{lemma} \label{l:L1l2-w-L1} R is $L^1(\ell^2))\to w-L^1$, i.e.,
\[  |\{x: |\sum_{j\in\Z}\psi_j(H)f_j(x)|>\al \}|\le C \al^{-1}\Vert 
(\sum_j |f_j|^2)^{1/2} \Vert_1
\]
\end{lemma}

\begin{proof} Let $\{f_j\} \in L^1(\ell^2)$, $\al >0 $. Let $ F(x)=
(\sum_{j=-\infty}^\infty |f_j(x)|^2 )^{1/2}$. By the C-Z decomposition for $F
\in L^1$ there exists a collection of disjoint open cubes $\{I_k\}$
such that   
\begin{align*}
i)&\quad |F(x)|\le \al, \;\;\textrm{a.e.} \;x \in \R^n\backslash \cup_k I_k\\
ii)& \quad\lam\le |I_k|^{-1}\int_{I_k } |F(x)|dx \le 2^n\al, \quad \forall k.
\end{align*}

Define 
$$ 
g_j(x)= \left\{
\begin{array}{ll}
|I_k|^{-1}\int_{I_k} f_j dx&  x\in I_k\\
f_j(x)                          & otherwise. 
\end{array}\right.  
$$
and 
$
b_j(x)= f_j(x)-g_j(x)$. Then $|g_j(x)|\le C\al$.

Let $x\in I_k$. Minkowski's inequality gives 
\begin{align*}
(\sum_{j=-\infty}^\infty& |g_j(x)|^2)^{1/2}= 
 (\sum_{j=-\infty}^\infty (|I_k|^{-1}\int_{I_k} |f_j| dx)^2)^{1/2}\\
\le& |I_k|^{-1}\int_{I_k} (\sum_{j=-\infty}^\infty  |f_j|^2 )^{1/2}dx
\le 2^n\al.
\end{align*}

Thus $\int \sum_j |g_j|^2\le c\al\Vert F\Vert_1$
and by Lemma \ref{l:Q-R-L2}
\begin{align*}
\vert&\{x: \sum_j |R_j g_j(x) |^2 )^{1/2} >\al/2 \}\vert 
\le C\al^{-2}
\Vert \{R_j g_j\}\Vert^2_{L^2(\ell^2)}\\
\le&
c\al^{-2} \Vert \{ g_j\}\Vert^2_{L^2(\ell^2)}
\le c \al^{-1} \Vert \{f_j\}\Vert_{L^1(\ell^2)}.
\end{align*}
(use the system $\{\phi_j\}_{-\infty}^\infty$)

It remains to estimate for $\{b_j(x)\}$.  Let $2^{j_k}\sim \ell(I_k)^2$. 


\begin{align*}
\int_{\R^n\setminus I_k^*}& (\sum_{j\in\Z} |R_j(1-\Phi_{j_k}(H)) b_{j,k}(x)|^2)^{1/2}dx\\
=& \int_{\R^n\setminus I_k^*} 
(\sum_{j\sim j_k} |\int_{I_k}\psi_j(H)(1-\Phi_{j_k}(H)(x,y)b_{j,k}(y)dy|^2)^{1/2}dx \;\\
\le& \int_{\R^n\setminus I_k^*}
\int_{I_k}(\sum_j |\psi_j(H)(1-\Phi_{j_k}(H)(x,y)b_{j,k}(y)|^2)^{1/2}dydx\quad \text{Minkowski} \\
\le& \int_{\R^n\setminus I_k^*}
\int_{I_k}\sup_j |\psi_j(H)(x,y)(1-\Phi_{j_k}(H)(x,y) | (\sum_j|b_{j,k}(y)|^2)^{1/2}dydx\quad \\
\le& \int_{I_k} (\sum_j |b_{j,k}(y)|^2)^{1/2} dy 
\int_{\R^n\setminus I_k^*}\sum_{j\ge j_k} |\psi_j(H)(1-\Phi_{j_k}(H)(x,y)|dx \\
\le& C \int_{I_k} (\sum_j |b_{j,k}(y)|^2)^{1/2}dy \le
C' \int_{I_k} (\sum_j |f_j(y)|^2)^{1/2}dy.
\end{align*}
where we used 
\[\int_{|x-y|>\ell(I_k)}|\psi_j(H)(x,y) |dx
\le (2^{j/2} \ell(I_k))^{-1} \qquad (N=n+1)
\]
by condition (\ref{e:int-phi-dec}). 

Hence, 
\begin{align*}
\vert&\{x\in\R^n\setminus \cup_kI_k^*: 
(\sum_j |R_j\sum_k(1-\Phi_{j_k}(H)) b_{j,k}(x) |^2)^{1/2} >\al/2 \}\vert\\
\le& 2\al^{-1}
\sum_k \int_{\R^n\setminus I_k^*} (\sum_j |R_j(1-\Phi_{j_k}(H))b_{j,k}(x)|^2)^{1/2} dx\\
\le& C \al^{-1} \int (\sum_j |f_j(y)|^2)^{1/2}dy,
\end{align*}
where $b_j=\sum_k b_{j,k}$ 
(convergence in $L^1\cap L^2$ so 
$T_j b_j(x)=  \sum_k T_jb_{j,k}$ in $L^2$)
and we used  
\[  (\sum_j  (\sum_k |T_jb_{j,k}(x)|)^2)^{1/2}
\le \sum_k(\sum_j |T_j b_{j,k}(x)|^2)^{1/2}\]
by Minkowski inequality.

To estimate $\vert\{x\notin\cup_kI_k^*:  (\sum_j |R_j\sum_k\Phi_{j_k}(H)) b_{j,k}(x) |^2)^{1/2} >\al/2 \}\vert$, enough 
\begin{equation}\label{e:sum-L2}
\sum_j\int_{\R^n}  | \psi_j(H)\sum_k \Phi_{j_k}(H) b_{j,k}(x)|^2 dx\le C \al  \Vert F\Vert_1. 
\end{equation}
Since $R=\{\psi_j(H)\}_{j\in\Z}$   
is uniformly bounded $L^2\to L^2$ or equivalently,
$ L^2(\ell^2)\to L^2(\ell^2)$, we only need show
\begin{equation}\label{e:sum-L2}
\sum_j\int_{\R^n} |\sum_k \Phi_{j_k}(H) b_{j,k}(x)|^2 dx\le C \al  \Vert F\Vert_1. 
\end{equation}
 Let $h=\{h_j\}\in L^2(\ell^2)$.
$\Vert h\Vert_{L^2(\ell^2)}\le 1$.  
\begin{align*}
&\sum_{j\in\Z} \la  \sum_k \Phi_{j_k}(H)  b_{j,k}(x)  , h_j\ra\\
\le& C \sum_j\sum_k \Vert b_{j,k}\Vert_1 |I_k|^{-1}\int_{z\in I_k}dz\int_ud\mu(u)\int_x \rho_{t}(x-z-u)
 |h_j(x)| dx\\ 
\le& C \sum_j \sum_k |I_k|^{-1} \Vert b_{j,k}\Vert_1\int_z\chi_{I_k}(z)  (\mu*M h_j)(z) dz\\
\le&C(\sum_j\Vert\sum_k |I_k|^{-1}\Vert b_{j,k}\Vert_1\chi_{I_k}\Vert_2^2)^{1/2}
\Vert M h_j\Vert_{L^2(\ell^2)}\\
\le & C(\sum_j  \sum_k  |I_k|^{-1} \Vert b_{j,k}\Vert_1^2 )^{1/2} \Vert h_j\Vert_{L^2(\ell^2)}\\
\le & C   [\sum_k |I_k|^{-1}\sum_j (\int_{I_k} |b_j(y)| dy )^2]^{1/2}\\ 
\le & C [\sum_k |I_k|^{-1} \big(\int_{I_k} (\sum_j |b_j |^2)^{1/2} dy \big)^2]^{1/2}
\le C  (\al\Vert F\Vert_1)^{1/2}. 
\end{align*}
\end{proof}
where we have used the Fefferman-Stein inequality: if $1<p<\infty$
\[   \Vert (\sum_j (M f_j)^2)^{1/2}\Vert_p\le \Vert (\sum_j |f_j|^2)^{1/2}\Vert_p \;.
\]

\nd
\rk In the above proof of two lemmas we can replace $\{\phi_j\}_{j=-\infty}^\infty$
with $\{\Phi, \phi_j\}_{j=1}^\infty$ to obtain the inhomogeneous result.
 The homogeneous result is necessary and useful for
obtaining Strichartz estimates for  wave equation.   

\section{P\"oschl-Teller potential}\label{S2}
For $H=-d^2/dx^2+V_\nu$ solve the Helmholtz equation 
\begin{equation}\label{e-evp}
He(x,z)=z^2 e(x,z).
\end{equation}
Under suitable asymptotic condition the solution also solves 
 the Lippman-Schwinger equation 
\begin{equation}\label{e-ls}
e (x,k)=e^{ikx}+\frac{1}{2i k}\int e^{i k\, |x-y|}V(y)e (y,k)\, dy\, ,
\end{equation}
that is 
\begin{align*}
 e (x,k)=&e^{ikx} - R_0(k^2+i0)Ve (\cdot,k) (x)\\
=& \sum_{n=0}^\infty (-R_0(k^2+i0)V)^n e^{ikx}\\
=& (I+R_0(k^2+i0)V)^{-1} e^{ikx}\;. \end{align*} 
Alternatively we can also using the above equations to write 
\begin{align*}
e(x,k)=& e^{ikx}-R_V(k^2)V e^{ikx} \,.
\end{align*}
where the free resolvent has the kernel 
\[  R_0(k^2+i0)(x,y)= -\frac{1}{2ik}e^{ik |x-y|} \qquad k\in \R\setminus \{0\} . 
\]

From \cite{OZ} we know that the continuous spectrum is  $\sigma_c=[0,\infty)$, and
the point spectrum is $\sigma_p={-1,-4,\cdots, -n^2}$. 
Bound states are 
Schwartz functions that are bounded by $c \, e^{-|x|}$. 








 We obtained in \cite{OZ} the following formula for the continuum eigenfunctions.
\begin{proposition}\label{th-exk}
Let $k\in \R\setminus \{0\}$. Then
$$e_n(x,k)=(\mathrm{sign}(k))^{n}\left(\prod_{j=1}^n\frac{1}{j+i | k|}
\right)\,P_n(x,k)e^{ikx},$$
where $P_n(x,k)=p_n(\th x,ik)$ is defined by the recursion formula
$$
p_n (\th x, ik)=
\frac{d}{dx}\big(p_{n-1}(\th x, ik)\big)+( ik-n\th x) p_{n-1}(\th x, ik)\, .
$$
In particular, the function
$$\R\times (\R\setminus \{0\})\ni (x,k)\mapsto e_n(x,k)\in\C$$
is analytic with $e_n(x,-k)=e_n(-x,k)$. 
Moreover, the function
$$(x,y,k)\mapsto e_n(x,k)\overline{e_n(y,k)}
=\left(\prod_{j=1}^n\frac{1}{j^2+ k^2}
\right)\, P_n(x,k)P_n(y,-k) e^{ik(x-y)}$$
is real analytic on $\R^3$.
\end{proposition}

\subsection{Resonance}   The Wronskian can be computed by using Jost functions    
\[
W(z)=-2(-1)^n ik \prod_{\iota=1}^n\frac{\iota+ik}
{\iota-ik} \,. \]

The remaining of this section is devoted to proving that conditions (a) and (b) in Theorem 
\ref{th:p-t-m(H)} are true for 
$H_\nu=H_0+V_\nu$\,. 


\subsection{Weighted $L^2$ for $m(H)(x,y)$}\label{ss:wei-L2-m(H)}  

Let $K_j(x,y)$ denote the kernel of $(m\phi_j)(H)E_{ac}$.  Let $\lam=2^{-j/2}$.  
\begin{lemma}
\begin{align}
&\Vert (x-y) K_j(\cdot,y)\Vert_2 \le \lam^{1/2}   \quad \forall j  \label{e:w-L2}\\  
&\Vert (x-y)\partial_yK_j(\cdot,y)\Vert_2 \le \lam^{-1/2} \quad  j\ge 0 \label{e:w-der-L2}    
\end{align}
\end{lemma}


\begin{proof} Proof of (\ref{e:w-L2}).   
\begin{align*}
&2\pi i(x-y) K_j(\cdot,y) \\
=&  i(x-y) \int m_j(k^2) e(x,k) \overline{e(y,k)}dk  \\
=&  \int_{|k|\sim\lam^{-1}}m_j(k^2)\left(\prod_{j=1}^n\frac{1}{j^2+ k^2}
\right)\, P_n(x,k)P_n(y,-k) \partial_k (e^{ik(x-y)}) dk\\
=&  -\int_{|k|\sim\lam^{-1}} \partial_k[m_j(k^2)\left(\prod_{j=1}^n\frac{1}{j^2+ k^2}
\right)\, P_n(x,k)P_n(y,-k)] e^{ik(x-y)} dk, 
\end{align*}
which can be written as finite sums of  
\[  (\th x)^\ell(\th y)^k  \big[ (m_j(k^2))' (\prod_{j=1}^n (j^2+k^2))^{-1} r_{2n}(k)\big]^\vee(x-y)
\]
\[(\th x)^\ell(\th y)^k  \big[ m_j(k^2) (\prod_{j=1}^n (j^2+k^2))^{-1} r_{2n-1}(k)\big]^\vee(x-y)
\]
and
\[ (\th x)^\ell(\th y)^k  \big[ m_j(k^2) \frac{2k}{\iota^2+k^2}(\prod_{j=1}^n (j^2+k^2))^{-1} 
q_{2n}(k)\big]^\vee(x-y)
\]
$0\le \ell, k,\iota\le n$, $0\le i\le 2n$, $r_i,q_i$ are polynomials of degree $i$.

Plancherel formula for Fourier transform gives
\begin{align*}
\Vert   (x-y)K_j(\cdot,y)\Vert_2
= O(\lam^{1/2})= O( 2^{-j/4})\qquad \forall j.
\end{align*}
using 
\[\begin{cases} 
(m_j(k^2))^{(i)}= O(\frac{1}{k^i})\quad i=0,1\\
 (\prod_{j=1}^n (j^2+k^2))^{-1} r_i(k)=O(1/\la k\ra)\\
 \frac{2k}{\iota^2+k^2}= O(1/\la k\ra)
\end{cases}
\]

\nd
Proof of (\ref{e:w-der-L2}).
\begin{align*}
&2\pi i(x-y) \partial_y K_j(\cdot,y) \\
=&  i(x-y) \int m_j(k^2) e(x,k) \overline{\partial_ye(y,k)}dk  \\
=& \int_{|k|\sim\lam^{-1}}(-ik)m_j(k^2)\left(\prod_{j=1}^n\frac{1}{j^2+ k^2}
\right)\, P_n(x,k)P_n(y,-k) \partial_k (e^{ik(x-y)}) dk\\
=&  \int_{|k|\sim\lam^{-1}}m_j(k^2)\left(\prod_{j=1}^n\frac{1}{j^2+ k^2}
\right)\, P_n(x,k)\partial_y(P_n(y,-k)) \partial_k (e^{ik(x-y)}) dk\\
:= I_1+I_2
\end{align*}

\begin{align*}
I_1=&  -\int_{|k|\sim\lam^{-1}} \partial_k[-ikm_j(k^2)\left(\prod_{j=1}^n\frac{1}{j^2+ k^2}
\right)\, P_n(x,k)P_n(y,-k)] e^{ik(x-y)} dk, 
\end{align*}

A similar argument as proving (\ref{e:w-L2}) yields
\begin{align*}
\Vert   (x-y)K'_{1,j}(\cdot,y)\Vert_2
= O(\lam^{-1/2})= O( 2^{j/4})\qquad \forall j.
\end{align*}
using 
\[\begin{cases} 
(m_j(k^2))^{(i)}= O(\frac{1}{k^i})\quad i=0,1\\
k= O(k)\\
 (\prod_{j=1}^n (j^2+k^2))^{-1} r_i(k)=O(1/\la k\ra)\\
 \frac{2k}{\iota^2+k^2}= O(1/\la k\ra)
\end{cases}
\]

\begin{align*}
I_2=&  -\int_{|k|\sim\lam^{-1}} \partial_k[m_j(k^2)\left(\prod_{j=1}^n\frac{1}{j^2+ k^2}
\right)\, P_n(x,k)\partial_y(P_n(y,-k))] e^{ik(x-y)} dk\\
&=  -\sech^2 y\int_{|k|\sim\lam^{-1}} \partial_k[m_j(k^2)\left(\prod_{j=1}^n\frac{1}{j^2+ k^2}
\right)\, P_n(x,k) (\partial_yp_n)(\th y,-k) ] e^{ik(x-y)} dk. 
\end{align*}
where we find if $k\sim \lam^{-1}$
\[
\partial_k [m_j(k^2)\left(\prod_{j=1}^n\frac{1}{j^2+ k^2}
\right)\, P_n(x,k) (\partial_yp_n)(\th y,-k)  ]
=\begin{cases}
O(k^{-2})=O(1) & |k|\to \infty\\
O(k^{-1})& |k|\to 0
\end{cases}\]

Plancherel formula for Fourier transform gives
\begin{align*}
\Vert   (x-y)K'_{2,j}(\cdot,y)\Vert_2
=\sech^2y\begin{cases}
 O(\lam^{3/2})=O(\lam^{-1/2})= O( 2^{j/4})&  j\in \N_0\\
O(\lam^{1/2})=O(2^{-j/4}) & j <0
\end{cases}
\end{align*}
\end{proof}

\begin{lemma}\label{l:j-ker-L2}     Let $j\in\Z$.
\begin{align*} 
&\Vert K_j(\cdot,y)\Vert_2 \le \lam^{-1/2}= 2^{j/4} 
\quad \forall y\label{e:kerj-L2}  \\
&\Vert K_j(\cdot,y)\Vert_\infty \lesssim  2^{j/2}\,.
\end{align*}

\end{lemma}


\begin{proof} $K_j(x,y)= (m\phi_j)(H)E_{ac}(x,y)$. 
\begin{align*}
&2\pi K_j(\cdot,y) =  \int m_j(k^2) e(x,k) \overline{e(y,k)}dk  \\
=&  \int_{|k|\sim\lam^{-1}}m_j(k^2)\left(\prod_{j=1}^n\frac{1}{j^2+ k^2}
\right)\, P_n(x,k)P_n(y,-k)  (e^{ik(x-y)}) dk
\end{align*}
which can be written as finite sums of  
\[  (\th x)^\ell(\th y)^k  \big[ (m_j(k^2)) (\prod_{j=1}^n (j^2+k^2))^{-1} r_{2n}(k)\big]^\vee(x-y)
\]
$\therefore$
\begin{align*}
&\Vert m_j(H)(x,y)\Vert_2 \le
\Vert m_j(k^2) (\prod_{j=1}^n (j^2+k^2))^{-1} r_{2n}(k)\Vert_{L^2_k}\\
\sim& \Vert m_j(k^2)\Vert_{L^2_k} 
\sim 2^{j/4}(\int_{|k|\le 1} |\phi(k^2)|^2dk )^{1/2}, \quad \forall j
\end{align*}
(if $m=1$, $\phi_j\subset [-2^{j},2^{j}]$), 
using
\[ \begin{cases} 
m_j(k^2)=O(1)\\
 r_{2n}(k)\prod_{j=1}^n (j^2+k^2)^{-1} =O(1)
\end{cases}
\]
\end{proof}

\begin{lemma}\label{l:j-der-ker-2-infty}  
\begin{align*} 
&\Vert(x-y)\partial_y K_j(\cdot,y)\Vert_2 \lesssim 
\begin{cases}
 2^{j/4}+\sech^2y 2^{-j/4} & j\to-\infty\\
2^{j/4} & j\to \infty
\end{cases}\\
&\Vert (x-y)\partial_yK_j(\cdot,y)\Vert_\infty \lesssim  
2^{j/2}+\sech^2y\begin{cases}
O(1) & j\to-\infty\\
2^{-j/2} & j\to \infty
\end{cases}
\end{align*}
\end{lemma}

\nd
\rk This means for $j\to -\infty$, 
$\Vert (x-y)\partial_y K_j(\cdot,y)\Vert_r\sim 
\sech^2y 2^{-j/(2r)}$, $r\in [2,\infty]$, which does not
seem to help establish the H\"ormander integral condition 
even if using $r$-norm instead of $2$-norm. 

\begin{proof} For $2$-norm, it is proved before.  For $r=\infty$,
\begin{align*}
&i 2\pi(x-y)\partial_yK_j(\cdot,y)\\
=& -\int \partial_k[km_j(k^2)\left(\prod_{j=1}^n\frac{1}{j^2+ k^2}
\right)\, p_n(\th x,k)p_n(\th y,-k)] e^{ik(x-y)} dk\\
+&(-i)\sech^2y\int_{|k|\sim\lam^{-1}} \partial_k[m_j(k^2)\left(\prod_{j=1}^n\frac{1}{j^2+ k^2}
\right)\, p_n(\th x,k)(\partial_yp_n)(\th y,-k)] e^{ik(x-y)} dk\\
\sim& \int_{|k|\sim\lam^{-1}}O(1)dk +\sech^2y 
 \int_{|k|\sim\lam^{-1}}\begin{cases}
O(1/k)& k\to 0\\ 
O(1/k^2)=O(1)&k\to\infty
\end{cases} dk
\end{align*}
($\lam=2^{-j/2}$)\\
$\therefore$
\begin{align*}
\Vert(x-y)\partial_yK_j(\cdot,y)\Vert_\infty
\sim 2^{j/2}+\sech^2y\begin{cases}
O(1) & j\to-\infty\\
2^{-j/2} & j\to \infty.
\end{cases}
\end{align*}\end{proof}

\subsection{Weighted pointwise decay of the kernel} 
The problem for pointwise decay of $\Phi_j(H)(x,y)$ in higher energy 
 can be overcome by using an integral version of (\ref{e:phi-dec}) with
a finite measure \footnote{ 
We modify Hebisch method when the kernel is rough
(and slowly decaying), not having Lipschitz smoothness
as needed in the H\"ormander method}. 

\begin{lemma}  (Hebisch-Zheng) \label{l:zheng-Harnack} 
Let $\Psi\in C^\infty_0$ be supported in $ [-1,1] $ and 
 let $I$ be any cube in $\R^n$ with length $\ell(I)$.    
 Then for all $x\in \R^n$ 
 and $y\in I$ with $\ell(I)=2^{-j/2}$ we have

a) \[ |\Phi_j(H)(x,y)|
 \le c \int_{u\in\R^n} \frac{2^{jn/2}}{(1+2^{j/2}| x-y-u|)^{n+\eps}} d\mu(u)\]
$d\mu(u)= \de(u)+\la u\ra^m e^{-c|u|} du$, some $m\ge 0$. 

Hence 
b)
\[ \sup_{y\in I} |\Phi_j(H)(x,y)|
 \le \frac{c}{|I|} \int_{z\in I}\int_{u\in\R^n}   
 \frac{2^{jn/2}}{(1+2^{j/2}| x-z-u|)^{n+\eps}} d\mu(u) dz\,.
\]
\end{lemma}

\begin{proof} Let $\Psi(x)=\Phi(x^2)$,  $\lam=2^{-j/2}$.  According to the 
formula for the kernel in preceding subsection 
\begin{align*}
 &\Phi_{j}(H)(x,y)\\
 =& \big[ \Psi (\lam k)\big]^\vee(x-y)
  +\big[ \Psi (\lam k)  
\sum_{\stackrel{\iota=1,\dots,n}{\mu_\iota=0,\dots,N} } 
\frac{a_{\mu_\iota,\iota} +b_{\mu_\iota,\iota} k}{(\iota^2+k^2)^{\mu_\iota+1} }  \big]^\vee(x-y)\\
=& \lam^{-1} {\Psi }^\vee(\lam^{-1} (x-y))\\
+&c\lam^{-1} \int_\R {\Psi }^\vee(\lam^{-1} (x-y-u)) |u|^m e^{- c|u|} du\,,\\
 \end{align*}
$m\in\N_0$ is an integral constant. Thus 

\begin{align*} 
&|\Phi_{j}(H)(x,y)|=  |\int_{|k|\lesssim 2^{j/2}} \Psi(2^{-j/2} k) e^{i(x-y)k}(1+\hat{a}(k)) dk |\\
=& \big(\lam^{-n} \Psi^\vee(\lam^{-1}\cdot) *(1+\hat{a} )^\vee\big) (x-y)\\
\lesssim& 
 \lam^{-n}\int_{\R^n}  (1+\lam^{-1} |x-y-u| )^{-n-\eps} d\mu(u)\;,
\end{align*}
where $y\in I$, $d\mu(u)=\de(u)+a(u)du$, $a(u)=|u|^m e^{-c |u|} $ .   
\end{proof}

\nd
\rk Observe that for $j>0$, $\Phi_j$ actually contains both high+low energy information
(if $0 \in \supp\,\Phi$).  This is the most technically difficult part.   Fortunately with Lemma 
\ref{l:zheng-Harnack} 
 we can control it by maximal function. 

\subsection{Kernel decay for a positive potential}
The case is simpler when $V$ is nonnegative. 
It can be shown \cite{Z06a, Da89} that  (\ref{e:phi-dec}) 
is a weaker 
assumption than the heat kernel estimate 
 \[  0\le   e^{-tH}(x,y)\le c_n t^{-n/2} e^{-c d(x,y)^2/t} \qquad \forall t>0\,.\] 
Examples include $H$ being a uniform elliptic operator or its perturbation of order $0$.
i.e., a Schr\"odinger operator with $V=V_+-V_-$,  $V_-$ is small in Kato norm,  cf. \cite{DP05}.

In \cite{Ou06} it is shown if $\sigma(H)\subset [0,\infty)$ then
\[ e^{-tH}(x,y)\le C t^{-n/2} e^{-|x-y|^2/4t} (1+ \de t+ |x-y|^2/t)^{n/2}\quad \forall t>0
\]
$\de=\de(V_-)$ and $\de=0$ if $V_-=0$.

\begin{proposition} Let $H$ denote a selfadjoint operator on $(M,g)$ with dimension $n$. 
Suppose $e^{-tH}$ verifies the upper Gaussian bound
\[  e^{-tH}(x,y)\le C_n t^{-n/2} e^{-c d(x,y)^2/t}\quad \forall t>0.
\]
Then for each $\ell$
\[    |\phi_j( H)(x,y)|\le C_\ell 2^{jn/2}(1+ 2^{j/2}d(x,y))^{-\ell}  \quad \forall j\in \Z\,. 
\]
\end{proposition}


\begin{thebibliography}{99}

\bibitem{Ag} S. Agmon, Spectral properties of Schr\"odinger operators
and scattering theory, {\em Annali Scuola Norm. Sup. de Pisa} {\bf 2}
(1975), 151--218.


\bibitem{A}
G.~Alexopoulos, Spectral multipliers on Lie groups of polynomial
growth, Proc. A.M.S. {\bf 120} (1994) no.3, 973-979.




\bibitem{BZ05} J.J.~Benedetto, S.~Zheng, Besov spaces for the 
Schr\"odinger operator with barrier potential (submitted).
{\em http://lanl.arXiv.org/math.CA/0411348}, (2005).


\bibitem{BL} J.~Bergh, J.~L\"ofstr\"om, 
{\em Interpolation Spaces}, Springer-Verlag, 
1976.

\bibitem{B99}
H.R.~Beyer, On the completeness of the quasinormal modes of the 
P\"oschl-Teller potential,
 {\em Comm. Math. Phys}. {\bf 204} (1999),  no. 2, 397-423. 
                                                                         

 
 







\bibitem{CS88}
F.M. Christ, C.D. Sogge, The weak type $L^1$ convergence of eigenfunction expansions for pseudodifferential operators, Invent. Math. 94 (2) (1988) 421--453.







\bibitem{CGM02}
M.~Cowling, A.~Sikora, A spectral multiplier theorem for a sublaplacian
on $SU(2)$,  Math. Z. 238 (2001), no. 1, 1--36.  


\bibitem
{DP05}
P. D'Ancona, V. Pierfelice, 
On the wave equation with a large rough potential. {\em J. F. A}. 
http://arXiv.org/math.AP/0310199.


\bibitem{Da89}
E.B. Davies, Heat Kernels and Spectral Theory, Cambridge University Press, Cambridge, 1989. 






\bibitem{DT79}
P.~Deift, E.~Trubowitz,  Inverse scattering on the line. \emph{Comm. Pure Appl. Math}. vol. XXXII.
(1979). 
121-251.


\bibitem{DM99} X.T.~Duong, A.~McIntosh, 
Singular integral operators with non-smooth kernels on irregular domains,
Rev. Mat. Iberoamericana,  {\bf 15} (1999), no.2: 233-265.

\bibitem
{DOS02} X.T.~Duong, E.M.~Ouhabaz and A.~Sikora, Plancherel type
estimates and sharp spectral multipliers, {\em J. Funct. Anal}. (2002). 

\bibitem{DY05}
 X.~T.~Duong,  L.~Yan, 
Duality of Hardy and BMO spaces associated with operators with heat kernel bounds, 
 J. Amer. Math. Soc. 18 (2005), 943-973. 



\bibitem{D99}
J.~Dziuba\'nsk,  Spectral multiplier theorem for $H^1$ spaces
associated with some Schr\"odinger operators, {\em Proc. A.M.S.} 
{\bf 127} (1999), no. 12, 3605-3613.


\bibitem{D01} \bysame, 
A spectral multiplier theorem for $H\sp 1$ spaces associated with Schr\"odinger operators with potentials satisfying a reverse H\"older inequality. 
{\em Illinois J. Math}. 45 (2001), no. 4, 1301--1313. 







\bibitem
{E96} 
 J.~Epperson, Hermite multipliers and pseudo-multipliers, {\em
Proc. Amer. Math. Soc.} {\bf 124} (1996), no. 7, 2061--2068.

\bibitem
{E97} ---, Hermite and Laguerre wave packet expansions. {\em Studia
Math.} {\bf 126} (1997), no. 3, 199--217.





\bibitem{Fe70}
C.~Fefferman, Inequalities for strongly singular convolution operators,
Acta Math. {\bf 124} (1970), 9-36.


\bibitem
{Flu74}
S.~Fl\"ugge, {\em Practical Quantum Mechanics}, Springer-Verlag, 1974.

\bibitem
{FJW} M. Frazier, B. Jawerth, G.Weiss, {\em Littlewood-Paley Theory
and the Study of Function Spaces},  Conference Board of the
Math. Sci. {\bf 79}, 1991.









\bibitem
{GH98} C.-A.~Guerin, M.~Holschneider, Time-dependent scattering
on fractal measures,  {\em J. Math. Physics} {\bf 39}(8), 1998.

\bibitem
{He90a} 
W.~Hebisch: A multiplier theorem for Schr\"odinger operators.
\textit{Colloq. Math}, {\bf 60/61} (1990),  no. 2, 659-664.

\bibitem{He90b} \bysame, 
 Almost everywhere summability of eigenfunction expansions associated to elliptic 
operators. \textit{Studia Math} {\bf 96} (1990), no. 3, 263--275.

\bibitem{He95} \bysame, Functional calculus for slowly decaying
kernels, preprint, 1995.\\ 
http://www.math.uni.wroc.pl/\symbol{126}hebisch









\bibitem
{JN94} A.~Jensen, S.~Nakamura, Mapping properties of functions
of Schr\"odinger operators between $L^p$ spaces and Besov spaces, 
in {\em Spectral and Scattering Theory and Applications}, Advanced 
Studies in Pure Math. {\bf 23}, 1994. 










\bibitem
{KMS} S.~Klainerman, M.~Machedon and J.~Stalker, Decay of
solutions to the wave equation on a spherically symmetric static 
background, {\em preprint}.




\bibitem
{Lam80}
G.L.~Lamb, Jr., {\em Elements of Soliton Theory}, 
Pure $\&$ Applied Mathematics, Wiley-Interscience, 1980.




\bibitem{Mi65}
S.G. Mikhlin, Multidimensional Singular Integrals and Integral Equations, Pergamon Press, Oxford, 1965.



\bibitem{MS94}
D. M\"uller, E.M. Stein, On spectral multipliers for Heisenberg and 
related groups, J. Math. Pures Appl. (9) 73 (4) (1994) 413--440. 







\bibitem
{OZ} G.~\'Olafsson, S.~Zheng, Function spaces associated with 
Schr\"odinger operators: the P\"oschl-Teller potential.  accepted for publication 
in {\em Journal of Fourier Analysis and Applications}. 

\bibitem{OOPSZ}
G.~\'Olafsson, K.~Oskolkov,
S.~Zheng, Spectral multipliers for  
Schr\"odinger operators: II,    Preprint.



\bibitem{Ou06}
E.M.~Ouhabaz, Sharp Gaussian bounds and $L^p$-growth of semigroups associated with
elliptic and Schr\"odinger operators. {\em Proc. A.M.S}. 2006. 

\bibitem{Ou06b} \bysame, Analysis of Heat Equations on Domains, London Math. Soc. Monographs, Vol. 31. Princeton Univ. Press 2004. 









\bibitem
{Sch05a} 
W.~Schlag, Dispersive estimates for Schr\"odinger operators: A survey. 
{\em http://lanl.arXiv.org/math.AP/0501037}, (2005). 

\bibitem
{Sch05b} 
\bysame, A remark on Littlewood-Paley theory for the distorted 
Fourier transform,
{\em http://lanl.arXiv.org/math.AP/0508577}, (2005)







\bibitem{Si82}
B.~Simon, Schr\"odinger semigroups, Bull. Amer. Math. Soc. {\bf 7}
(1982) no.3, 447-526.








\bibitem{St93}
E.~Stein, Harmonic analysis, Real-variable methods, orthogonality, and 
oscillatory integrals, Princeton Univ. Press, 1993.






\bibitem{Tay} 
M.~Taylor, $L^p$-estimates on functions of the Laplace operator, 
{\em Duke Math. J}. {\bf 58} no. 3 (1989), 773Ð793





\bibitem
{Tr83} H.~Triebel, {\em Theory of Function Spaces},
Birkh\"{a}user Verlag, 1983.

\bibitem
{Tr92} ---, {\em Theory of Function Spaces II}, Monographs
Math. {\bf 84}, Birkh\"auser, Basel, 1992.



\bibitem{W99}
R. Weder,  The $W^{k,p}$-continuity of the Schr\"odinger wave operators on the line. 
Comm. Math. Phys. {\bf 208} (1999) 507-520. 

\bibitem{Y95}
Yajima, K., The $W^{k,p}$-continuity of wave operators for Schr\"odinger 
operators. J. Math. Soc. Japan 47, 551--581 (1995).



\bibitem{Y05} 
\bysame, Dispersive estimate for Schr\"odinger equations with threshold
 resonance and eigenvalue.  Comm. Math. Phys. 2005.

\bibitem
{Zha} Q.~Zhang, Global bounds of Schr\"odinger heat kernels
with negative potentials.  {\em J. Func. Anal.} {\bf 182} (2001), no.2, 344-370. 

 
\bibitem
{Z04} S.~Zheng, A representation formula related to Schr\"odinger 
operators.  {\em Anal.~Theo.~Appl.} {\bf 20} (2004), no.3. 
{\em http://lanl.arXiv.org/math.SP/0412314}.


\bibitem
{Z06a} \bysame,  Littlewood-Paley theorem for  Schr\"odinger operators. (sumitted) 
2006. 

\bibitem
{Z06b}  
\bysame, Spectral multipliers, function spaces and dispersive estimates for Schr\"odinger 
operators.  
{\em Preprint}. 




\end{thebibliography}
\end{document}